\input{aipcheck}

\documentclass[
    ,final            
  ]
  {aipproc}
  
  \usepackage[centertags]{amsmath}
\usepackage{latexsym}
\usepackage{amsfonts}
\usepackage{amssymb}
\usepackage{amsthm}
  
\newenvironment{prof}[1][Proof]{\textbf{#1.} }{\ \rule{0.5em}{0.5em}}

\newtheorem{theorem}{Theorem}[section]
\newtheorem{definition}{Definition}[section]

\layoutstyle{6x9}

\begin{document}

\title{Conservation laws for under determined systems of differential equations}


\author{Mahouton Norbert Hounkonnou}{ address={University of Abomey-Calavi, International Chair in
    Mathematical Physics and Applications (ICMPA-UNESCO
    Chair), 072 B.P.: 50, Cotonou,
    Republic of Benin} \\
  Email: norbert.hounkonnou@cipma.uac.bj\footnote{With copy to
    hounkonnou@yahoo.fr.}}

\author{Pascal Dkengne Sielenou}{ address={University of Abomey-Calavi, International Chair in
    Mathematical Physics and Applications (ICMPA-UNESCO
    Chair), 072 B.P.: 50, Cotonou,
    Republic of Benin} }

\begin{abstract}
This work extends the Ibragimov's  conservation theorem for partial differential equations  [{\it J. Math. Anal. Appl. 333 (2007 311-328}] 
to under determined systems of differential equations.
   The concepts of adjoint equation and formal Lagrangian for a
 system of differential equations whose the number of equations is equal to or
lower than the number of dependent variables are defined.
 It is proved that the system given by an equation and its  adjoint
is associated with a variational problem (with or without classical Lagrangian) and inherits all Lie-point and generalized symmetries from the original equation.
 Accordingly, a Noether theorem for conservation laws can be formulated.
\end{abstract}

\maketitle

\section{Introduction}

Conservation laws play a vital role in the study of partial differential equations (PDEs) in the search for  qualitative properties
 such as integrability, stability,  existence of global solutions  and  the linearizability
conditions \cite{rm, r1,r2,r3,r4}. Their usefulness has considerably increased since the work by Jacobi \cite{s1} 
 in 1884, who 
showed a connection between conserved quantities and symmetries of the
equations of a particle's motion in classical mechanics.
  Klein  \cite{s2}  has obtained similar result
 for the equations of the general relativity and predicted that a connection between conservation laws and symmetries could be
found for any differential equation obtained from a variational principle.
 Noether \cite{r8} has showed that  the conservation laws were associated with the invariance of variational integrals with respect to continuous transformation groups. She obtained the sufficient condition for existence of conservation laws.
 In 1921,  Bessel-Hagen  \cite{s3}   applied Noether's theorem with the so-called "divergence" condition to
the Maxwell equations and calculated their conservation laws.
 In 1951, Hill \cite{s4}  presented the explicit formula
 in terms of variations for conservation laws in the case of a first order Lagrangian.
  Ibragimov  \cite{s5}  proved the generalized version of Noether's theorem  and   conservation laws  related to the invariance of the extremal values of
variational integrals. He derived  the necessary and sufficient condition for the
existence of conservation laws and gave the explicit expressions
 in the case of a Lagrangian of any order.

Recently, some relevant results on conservation laws have emerged. See for instance the works by 
Wolf \cite{r12}, Anco \cite{r6, r7},
  Poole \cite{r16},
Ibragimov \cite{r13a, r13},
Naz \cite{r13b} and
     Khamitova  \cite{r13c}.
Despites all this progress, there remains an important question:
How to associate a conservation law with every infinitesimal generator of
symmetries of an arbitrary differential equation?
  Ibragimov \cite{r13} achieved this goal for any system of differential equations
where the number of equations is equal to the number of dependent variables.
 Our work extends Ibragimov's result to under determined system of differential equations.
\section{Notations, basic definitions and theorems}
Consider $X$, an $n$-dimensional independent variables space, and $U=\bigotimes_{j=1}^{m}U^j,$
 an $m$-dimensional dependent variables space.
 Let $x=(x^1,\cdots,x^n)\in X$ and $u=(u^1,\cdots,u^m)\in U$ with $u^j\in U^j.$
  We define  the Jet-space $U^{(s)}$ as
  $ U^{(s)}:=\bigotimes_{j=1}^{m}\left(\bigotimes_{l=0}^{s}U^j_{(l)}  \right),$
 where $U^j_{(l)}$ is the set of all $p_l\equiv \left(\begin{array}{c}{n+l-1}\\{l}  \end{array}\right)$ distinct $l$-th order
 partial derivatives of $u^j.$
 We denote by $u^j_{(k)}$ the $p_k$-tuple of all $k$-order  derivatives of $u^j$.
 An element   $u^{(s)}$ in the Jet-space $U^{(s)}$  is the $m(1+p_1+p_2+\cdots+p_s)=m \left(\begin{array}{c}{n+s}\\{s}  \end{array}\right)$-tuple defined by
$u^{(s)}=\left(u^1_{(0)},u^1_{(1)},\cdots,u^1_{(s)}, u^2_{(0)},u^2_{(1)},\cdots,u^2_{(s)},\cdots,u^m_{(0)},u^m_{(1)},\cdots,u^m_{(s)}  \right).$

A variational problem consists in finding extrema of a functional $\mathfrak{L}$ defined by
\begin{equation}\label{eq6}
\mathfrak{L}[u]=\int_{\Omega}L(x,u^{(s)})dx,
\end{equation}
where $\Omega$ is a connected open subset of $X$ and $L$ defined on $X\times U^{(s)}$ is an  $s$-order differential function
 called the Lagrangian of the variational problem $\mathfrak{L}.$
 In general, a functional is a mapping that assigns to each element in some
function space a real number, and a variational problem amounts to searching for functions which are an extremum (minimum, maximum, or saddle
points) of a given functional.
\begin{theorem}
Let $u$ be an extremum of $\mathfrak{L},$ then $u$ satisfies the Euler-Lagrange equations
\begin{equation}\label{eq7}
 \frac{\delta}{\delta u^j} L(x,u^{(s)})  =0,\quad j=1,\cdots,m.
\end{equation}
\end{theorem}
\begin{theorem}\textbf{(Noether's theorem)}\\
 Let $G$ be a one parameter variational symmetry group for the functional $\mathfrak{L}[u]=\int_{\Omega}L(x,u^{(s)})dx,$
and
$V=\sum_{i=1}^n\xi^i(x,u)\frac{\partial}{\partial x^i}+\sum_{j=1}^m\phi^j(x,u)\frac{\partial}{\partial u^j}$
 be an infinitesimal generator  of $G,$ i.e. $V$ satisfies the variational infinitesimal invariance criterion
 \begin{eqnarray}\label{eqa0}
 Pr^{(s)}V(L)+L\,div\,\xi=div\, B
 \end{eqnarray}
 for some vector $B=(B^1,\cdots,B^n)$ of differential functions, where $Pr^{(s)} X$ is the $s-$order prolongation of $X$. 
Then the vector $C=(C^1(x,u^{(s_1)}),\cdots,C^n(x,u^{(s_n)}))$ defined by:
\begin{eqnarray}
  C^i&=&-B^i+\xi^iL
  +\sum_{j=1}^m\sum_{k_1=0}^{s^j_1}\cdots \sum_{k_i=0}^{s^j_i-1} \cdots\sum_{k_n=0}^{s^j_n} D^{k_1}_{x^1}\cdots D^{k_n}_{x^n}(W^j)\nonumber\\
& \times& \sum_{l_1=0}^{s^j_1-k_1}\cdots\sum_{l_n=0}^{s^j_n-k_n}(-D_{x^1})^{l_1}\cdots 
 (-D_{x^n})^{l_n}         \left(\frac{\partial L}{\partial u^j_{(k_1+l_1)x^1\cdots (k_i+l_i+1)x^i\cdots (k_n+l_n)x^n}}\right)\nonumber
\end{eqnarray}
where
$W^j=\phi^j-\sum_{i=1}^{n}\xi^iu^j_{x^i},\quad j=1,\cdots,\,m$
provides a conservation law for the Euler-Lagrange equations (\ref{eq7}):
 $  \frac{\delta}{\delta u^j} L(x,u^{(s)})  =0,\quad j=1,\cdots,\, m,$
i.e. obeys the equation
$div\,C\equiv D_{x^1}\,C^1+\cdots + D_{x^n}\,C^n=0$
for all solution of the system (\ref{eq7}).
Such a vector $C$ 
is called a \textbf{conserved vector} for the system (\ref{eq7}).
\end{theorem}
 \section{Main results for under determined systems of nonlinear partial differential equations}
\begin{definition}\textbf{(Adjoint equation)}

Consider the system
\begin{equation}\label{eqa1}
  F_{\alpha}\left( x,u^{(s)},\widetilde{u}^{(s)} \right)=0,\quad \alpha=1,\ldots,m,
  \end{equation}
 where $F_{\alpha}$ are differential functions having $n$ independent variables
 $x=(x^1,\ldots,x^n)$ and $m+\widetilde{m}$ dependent variables $u=(u^1,\ldots,u^m),$
 $\widetilde{u}=(\widetilde{u}^1,\ldots,\widetilde{u}^{\widetilde{m}}),$ $u=u(x),$ $\widetilde{u}=\widetilde{u}(x);$
 $u^{(s)}$ (resp. $\widetilde{u}^{(s)}$ ) is a vector encompassing dependent variable $u$ (resp. $\widetilde{u}$ ) and their
  derivatives up to order $s.$
 We introduce the differential functions
\begin{eqnarray}
  F^{\star}_{\alpha}&=& \frac{\delta}{\delta u^{\alpha}}\left[\sum_{\beta=1}^m v^{\beta}F_{\beta}
  +\left( \sum_{\widetilde{\beta}=1}^{\widetilde{m}}\widetilde{v}^{\widetilde{\beta}} \right)\left( \sum_{\nu=1}^m F_{\nu} \right)  \right]  \nonumber\\
 \widetilde{F}^{\star}_{\widetilde{\alpha}}&=& \frac{\delta}{\delta \widetilde{u}^{\widetilde{\alpha}}}\left[\sum_{\beta=1}^m v^{\beta}F_{\beta}
  +\left( \sum_{\widetilde{\beta}=1}^{\widetilde{m}}\widetilde{v}^{\widetilde{\beta}} \right)\left( \sum_{\nu=1}^m F_{\nu} \right)  \right] \nonumber,
 \end{eqnarray}
 where $v=(v^1,\ldots,v^m)$ and $\widetilde{v}=(\widetilde{v}^1,\ldots,\widetilde{v}^{\widetilde{m}})$ are new dependent variables,
 ($v=v(x),$ $\widetilde{v}=\widetilde{v}(x)$), also called nonlocal variables. Then, we
define the corresponding system of adjoint equations 
by
\begin{eqnarray}\label{eqa2}
  F^{\star}_{\alpha}\left( x,u^{(s)},\widetilde{u}^{(s)},v^{(s)},\widetilde{v}^{(s)} \right)&=& 0,\quad \alpha=1,\ldots,m\\
  \label{eqa3}
 \widetilde{F}^{\star}_{\widetilde{\alpha}}\left( x,u^{(s)},\widetilde{u}^{(s)},v^{(s)},\widetilde{v}^{(s)} \right)&=& 0, \quad \widetilde{\alpha}=1,\ldots,\widetilde{m}.\label{eqa02}
 \end{eqnarray}
 \end{definition}
 \begin{theorem}
Any system of PDEs  (\ref{eqa1}):
$
  F_{\alpha}\left( x,u^{(s)},\widetilde{u}^{(s)} \right)=0,\quad \alpha=1,\ldots,m,
  $
considered together with their adjoint equations (\ref{eqa2})-(\ref{eqa02}),
has a Lagrangian. Namely, the  Eqs (\ref{eqa1})-(\ref{eqa02}) with $2(m+\widetilde{m})$ unknowns
 are the Euler-Lagrange equations with Lagrangian
$$
 L\left( x,u^{(s)},\widetilde{u}^{(s)},v^{(s)},\widetilde{v}^{(s)} \right)=
 \sum_{\beta=1}^m v^{\beta}F_{\beta}
  +\left( \sum_{\widetilde{\beta}=1}^{\widetilde{m}}\widetilde{v}^{\widetilde{\beta}} \right)\left( \sum_{\nu=1}^m F_{\nu} \right).
 $$
\end{theorem}
\begin{prof} It is immediate from the definitions of Euler-Lagrange equations and adjoint equations.
\end{prof}
\begin{definition}
The system (\ref{eqa1})
is called self-adjoint if the substitution $(v,\widetilde{v})=(u,\widetilde{u})$ in its adjoint Eqs. (\ref{eqa2})-(\ref{eqa02})
gives, for some differential functions $\Gamma_{\alpha\nu}$ and $\widetilde{\Gamma}_{\widetilde{\alpha}\nu},\,$
 $ F^{\star}_{\alpha}= \sum_{\nu=1}^m \Gamma_{\alpha\nu}F_{\nu} ,\quad \alpha=1,\ldots,\,m \label{eqb},\,\,$ and
 $\widetilde{F}^{\star}_{\widetilde{\alpha}}= \sum_{\nu=1}^m \widetilde{\Gamma}_{\widetilde{\alpha}\nu}F_{\nu},
\quad \widetilde{\alpha}=1,\ldots,\,\widetilde{m}\label{eqc}$.
The system (\ref{eqa1}) is called quasi-self-adjoint if there exist two functions $h$ and $\widetilde{h}$ such 
that the same expansions of $F^{\star}_{\alpha}$ and $\widetilde{F}^{\star}_{\widetilde{\alpha}}$ in terms of
 $\Gamma_{\alpha\nu}$ and $\widetilde{\Gamma}_{\widetilde{\alpha}\nu}\,$
 hold upon the substitution $(v,\widetilde{v})=\left(h(u,\widetilde{u}), \widetilde{h}(u,\widetilde{u})\right).$
\end{definition}
Provided these statements, we can now provide the main results of this paper.
\begin{theorem}
 Consider the system (\ref{eqa1}). Then
Its adjoint Eqs. (\ref{eqa2})-(\ref{eqa02})
 inherits symmetries of equations (\ref{eqa1}).
Namely, if the system (\ref{eqa1}) admits an operator
$X=\sum_{i=1}^n\xi^i\frac{\partial}{\partial x^i}+\sum_{\beta=1}^m \eta_{\beta}\frac{\partial}{\partial u^{\beta}}
+\sum_{\widetilde{\beta}=1}^{\widetilde{m}} \widetilde{\eta}_{\widetilde{\beta}}\frac{\partial}{\partial \widetilde{u}^{\widetilde{\beta}}}$,
where $X$ is a generator of a point transformation group, i.e. $\xi^i=\xi^i(x,u,\widetilde{u}),$ 
$\eta_{\beta}=\eta_{\beta}(x,u,\widetilde{u}),$ $\widetilde{\eta}_{\widetilde{\beta}}=\widetilde{\eta}_{\widetilde{\beta}}(x,u,\widetilde{u})$
and
$
Pr^{(s)}X(F_{\alpha})=\sum_{\beta=1}^m\lambda_{\alpha\beta}(x,u^{(s)},\widetilde{u}^{(s)})\,F_{\beta},
$
then the equations (\ref{eqa2})-(\ref{eqa02}) have the generator of symmetries
$Y=X+\sum_{\beta=1}^m \eta^{\star}_{\beta}\frac{\partial}{\partial v^{\beta}}+\sum_{\widetilde{\beta}=
1}^{\widetilde{m}} \widetilde{\eta}^{\star}_{\widetilde{\beta}}\frac{\partial}{\partial \widetilde{v}^{\widetilde{\beta}}}$,
where
$\eta^{\star}_{\beta}+\sum_{\widetilde{\alpha}=1}^{\widetilde{m}} \widetilde{\eta}^{\star}_{\widetilde{\alpha}}=
-\left[ \sum_{\alpha=1}^m v^{\alpha}\lambda_{\alpha\beta}+ v^{\beta}\sum_{i=1}^n D_{x^i}(\xi^i)
\label{eqa9}
+
\left(\sum_{\widetilde{\alpha}=1}^{\widetilde{m}} \widetilde{v}^{\widetilde{\alpha}}\right)\sum_{\alpha=1}^m
 \lambda_{\alpha\beta}+\left(\sum_{\widetilde{\alpha}=1}^{\widetilde{m}} \widetilde{v}^{\widetilde{\alpha}}\right)\left(\sum_{i=1}^n D_{x^i}(\xi^i)\right)  \right]$
which is satisfied in particular for
$\eta^{\star}_{\beta}= -\left[ \sum_{\alpha=1}^m\left(v^{\alpha} +\sum_{\widetilde{\alpha}=1}^{\widetilde{m}} \widetilde{v}^{\widetilde{\alpha}} 
\right)\lambda_{\alpha\beta} +v^{\beta}\sum_{i=1}^n D_{x^i}(\xi^i) \right],\nonumber\\
 \widetilde{\eta}^{\star}_{\widetilde{\alpha}}= - \widetilde{v}^{\widetilde{\alpha}}\sum_{i=1}^n D_{x^i}(\xi^i).$
\end{theorem}
\begin{prof} Using the variational infinitesimal test and setting the coefficients of $F_{\beta}$ to 0 yield the results.
\end{prof}
\begin{theorem}
 Consider the system (\ref{eqa1})
Its adjoint Eqs. (\ref{eqa2})-(\ref{eqa02})
 inherits Lie-B\"{a}cklund operator of equations (\ref{eqa1}).
Namely, if the system (\ref{eqa1}) admits an operator
$X=\sum_{i=1}^n\xi^i\frac{\partial}{\partial x^i}+\sum_{\beta=1}^m \eta_{\beta}\frac{\partial}{\partial u^{\beta}}
+\sum_{\widetilde{\beta}=1}^{\widetilde{m}} \widetilde{\eta}_{\widetilde{\beta}}\frac{\partial}{\partial \widetilde{u}^{\widetilde{\beta}}},$
where $X$ is a Lie-B\"{a}cklund operator, i.e. $\xi^i=\xi^i(x,u^{(s_1)},\widetilde{u}^{(\widetilde{s_1})}),$ 
$\eta_{\beta}=\eta_{\beta}(x,u^{(s_2)},\widetilde{u}^{(\widetilde{s_2})}),$ $\widetilde{\eta}_{\widetilde{\beta}}=
\widetilde{\eta}_{\widetilde{\beta}}(x,u^{(s_3)},\widetilde{u}^{(\widetilde{s_3})})$
are any differential functions,
and
$
Pr\,X\left(  F_{\nu} \right)=\sum_{\mu=1}^m\mathcal{D}_{\nu\mu}(F_{\mu})
$
for some differential operators
$\mathcal{D}_{\nu\mu}= \lambda^0_{\nu\mu}+\sum_{i_1=1}^n\lambda^{i_1}_{\nu\mu}D_{x^{i_1}}+\sum_{i_1,i_2=1}^n\lambda^{i_1i_2}_{\nu\mu}D_{x^{i_1}}D_{x^{i_2}}
    + \sum_{i_1,i_2,i_3=1}^n\lambda^{i_1i_2i_3}_{\nu\mu}D_{x^{i_1}}D_{x^{i_2}}D_{x^{i_3}}+\cdots,$
then equations (\ref{eqa2})-(\ref{eqa02}) admit the Lie-B\"{a}cklund operator
$Y=X+\sum_{\beta=1}^m \eta^{\star}_{\beta}\frac{\partial}{\partial v^{\beta}}+\sum_{\widetilde{\beta}=1}^{\widetilde{m}} 
\widetilde{\eta}^{\star}_{\widetilde{\beta}}\frac{\partial}{\partial \widetilde{v}^{\widetilde{\beta}}},$
where
\begin{eqnarray}
\eta^{\star}_{\beta}+\sum_{\widetilde{\beta}=1}^{\widetilde{m}} \widetilde{\eta}^{\star}_{\widetilde{\beta}}&=&
-\left\{w^{\beta}\sum_{i_1=1}^n D_{x^{i_1}}\left(\xi^{i_1}\right)
  +\sum_{\mu=1}^m\left[w^{\mu}\lambda_{\mu\beta}^{0} -\sum_{i_1=1}^n D_{x^{i_1}}\left(w^{\mu}\lambda_{\mu\beta}^{i_1}\right)
 \right.\right.\nonumber\\
  &+&\sum_{i_1,i_2=1}^nD_{x^{i_1}}D_{x^{i_2}}\left(w^{\mu}\lambda_{\mu\beta}^{i_1i_2}\right)
 -\left.\left.\sum_{i_1,i_2,i_3=1}^nD_{x^{i_1}}D_{x^{i_2}}D_{x^{i_3}}\left(w^{\mu}\lambda_{\mu\beta}^{i_1i_2i_3}\right) +\cdots  \right]\right\}
\nonumber
\end{eqnarray}
with $w^{\alpha}=v^{\alpha}
  +\sum_{\widetilde{\beta}=1}^{\widetilde{m}}\widetilde{v}^{\widetilde{\beta}}.$
  \end{theorem}
  \begin{prof} Replacing $w^{\beta}=v^{\beta}
  +\sum_{\widetilde{\beta}=1}^{\widetilde{m}}\widetilde{v}^{\widetilde{\beta}}$ in the result of the computation of
$T\equiv \emph{\emph{Pr}} Y(L)+L \emph{\emph{Div}}(\xi)$  and setting the coefficients of $F_{\beta}$ to 0, then the resulting expression of  $T$ reduces to a divergence, what
 achieves the proof.
\end{prof}

  Finally, we arrive at the following  general formulation of the conservation theorem.
\begin{theorem}
Every infinitesimal generator (or Lie-B\"{a}cklund operator)
\begin{equation}
X=\sum_{i=1}^n\xi^i\frac{\partial}{\partial x^i}+\sum_{\beta=1}^m \eta_{\beta}\frac{\partial}{\partial u^{\beta}}+\sum_{\widetilde{\beta}=1}^{\widetilde{m}} \widetilde{\eta}_{\widetilde{\beta}}\frac{\partial}{\partial \widetilde{u}^{\widetilde{\beta}}}
\end{equation}
  of differential equations (\ref{eqa1})
provides a conservation law for the system of differential equations
 (\ref{eqa1}) and their adjoint  (\ref{eqa2})-(\ref{eqa02}).
\end{theorem}
 Relevant  nonlinear partial differential equations of mathematical physics such as
the three dimensional time-dependent Ginzburg-Landau (GL) equations \cite{r15, r14} are studied in this framework. Main results of these investigations will be in the core of forthcoming papers.
\begin{theacknowledgments}
  This work is partially supported by the ICTP through the
OEA-ICMPA-Prj-15. The ICMPA is in partnership with
the Daniel Iagolnitzer Foundation (DIF), France.
\end{theacknowledgments}



\bibliographystyle{aipproc}   


\IfFileExists{\jobname.bbl}{}
 {\typeout{}
  \typeout{******************************************}
  \typeout{** Please run "bibtex \jobname" to optain}
  \typeout{** the bibliography and then re-run LaTeX}
  \typeout{** twice to fix the references!}
  \typeout{******************************************}
  \typeout{}
 }


\end{document}